\documentclass[12pt]{amsart}
\usepackage{amsmath,amsfonts,amstext,amsthm,epsfig}

\newtheorem{Theorem}{Theorem}
\newtheorem{Prop}{Proposition}
\newtheorem{Def}{Definition}
\makeatletter
    \addtocounter{section}{0}
    
     \@addtoreset{equation}{section}
\makeatother

\def\N{\mathbb N}
\def\R{\mathbb R}

\begin{document}

\title[Backward uniqueness]
{On the backward uniqueness property for a class of parabolic operators}%
\author{Daniele Del Santo}\author{Martino Prizzi}

\address{Daniele Del Santo, Universit\`a di Trieste, Dipartimento di
Matematica e Informatica,
Via Valerio 12/1, 34127 Trieste, Italy}%
\email{delsanto@univ.trieste.it}%

\address{Martino Prizzi, Universit\`a di Trieste, Dipartimento di
Matematica e Informatica, Via Valerio 12/1, 34127 Trieste, Italy}%
\email{prizzi@dsm.univ.trieste.it}%

\date{\today}%
\begin{abstract}
We give sharp regularity conditions,
ensuring the backward uniquess property to a class of parabolic operators.
\end{abstract}
\maketitle

\section {Introduction, statements and remarks}

In this note we illustrate  some  new results concerning the backward uniqueness property for a class of
parabolic operators, whose coefficients are non-Lipschitz continuous in time.
Namely, we consider parabolic operators of the form
\begin{equation}\label{operator}
P:=\partial_t+{\hskip-8pt}\sum_{0\leq|\alpha|,|\beta|\leq m}{\hskip-8pt}
(-1)^{|\alpha|}\partial_x^{\alpha}(\rho_{\alpha\beta}(t,x)\partial_x^{\beta});
\end{equation}
here $m\in\N$, $(t,x)\in[0,T]\times\R^n$, and $\alpha$ and $\beta$ are $n$-multiindices with weights $|\alpha|$ and
$|\beta|\leq m$. We assume that
$\rho_{\alpha\beta}=\overline{\rho_{\beta\alpha}}$ for all $\alpha$'s and $\beta$'s (formal self-adjointness), that
$\rho_{\alpha\beta}$ is real when $|\alpha|=|\beta|=m$, and that there exists $c_0>0$ such that
$\sum_{|\alpha|=|\beta|=m}\rho_{\alpha\beta}(t,x)\xi^\alpha\xi^\beta\geq c_0|\xi|^{2m}$ for all $\xi\in\R^n$ (strong
ellipticity).

\par Given a
functional space
${\mathcal H}$, we say that
$P$ enjoys the {\it backward uniqueness property in ${\mathcal H}$} iff, whenever $u\in{\mathcal H}$ satisfies $Pu\equiv 0$
(in the sense of distributions) in $[0,T]\times\R^n$, and
$u(T,\cdot)\equiv 0$ in $\R^n$, then $u\equiv 0$ in $[0,T]\times\R^n$. Our aim is to find conditions on the coefficients
$\rho_{\alpha\beta}$'s, ensuring that $P$ enjoys the backward uniqueness property in some given functional space $\mathcal
H$. As a preliminary observation, we notice that in
\cite{Ty} Tychonoff constructed a function $u\in C^\infty (\R\times\R^n)$ satisfying
\begin{equation}
\begin{cases}
\partial_t u-\Delta u\equiv 0&\text{in $\R\times\R^n$}\\
u(0,\cdot)\equiv 0&\text{in $\R^n$}\end{cases}
\end{equation}
but $u\not\equiv 0$ in any open subset of $\R\times\R^n$. It follows that, whether $P$ enjoys the {\it backward uniqueness
property in ${\mathcal H}$} or not, depends first of all on the choice of $\mathcal H$. We are interested here in the
case
${\mathcal H}={\mathcal H}_1^m$, where
\begin{equation}
{\mathcal H}_1^m:=H^1([0,T], L^2(\R^n))\cap L^2([0,T], H^{2m}(\R^n)).
\end{equation}
The reason for this choice is essentially due to its historical
background, but other choices are possible as well. In \cite{LM}
Lions and Malgrange proved that $P$ enjoys the backward uniqueness
property in ${\mathcal H}_1^m$, provided the
$\rho_{\alpha\beta}$'s are sufficiently smooth with respect to $x$
and Lipschitz continuous with respect to $t$. They work in an
abstract Hilbert space setting and their proof is based on a
Carleman type estimate. The required smoothness of the
$\rho_{\alpha\beta}$'s with respect to $x$ is related to the
regularity theory for elliptic equations and is needed to let $P$
fall in the abstract Hilbert space setting. The required Lipschitz
continuity with respect to $t$ seems to be more intrinsically
connected with the backward uniqueness property. Infact, in the
same paper Lions and Malgrange raised the question, whether
Lipschitz continuity could be replaced by, say, simple continuity.
As a first step in this direction, in \cite{BT} Bardos and Tartar
proved that $P$ enjoys the backward uniqueness property in
${\mathcal H}_1^m$, provided the $\rho_{\alpha\beta}$'s are
absolutely continuous with respect to $t$. Their proof exploits a
sort of {\it logarithmic convexity property} satisfied by the norm
$\|u\|$ of any nontrivial solution of $Pu=0$. Later, in \cite{Gh}
Ghidaglia, by using the same technique, extended the results of
\cite{BT}, so as to cover also some classes of nonlinear parabolic
equations.  We stress that, in all the above mentioned results, it
is required that the $\rho_{\alpha\beta}$'s be differentiable with
respect to $t$, at least in a weak sense. The reason is that, at a
certain point, one needs to perform some {\it integration by
parts}. Although this latter seems to be just a technical
obstruction, the possibility of replacing Lipschitz continuity by
simple continuity was finally ruled out by Miller in \cite{mill}.
He exibited an example of an operator $P$ which does not enjoy the
backward uniqueness property in ${\mathcal H}_1^m$. The operator
is of second order in space and its coefficients are of class
$C^\infty$ with respect to $x$ and H\"older continuos of exponent
$1/6$ with respect to $t$.
Recently, in \cite{M}, Mandache improved the result of Miller
constructing a similar non-uniqueness example in which the
coefficients are  of class $C^\infty$ with respect to $x$ and
H\"older continuos of every exponent less than 1 with respect to
$t$. More precisely in the result of Mandache the regularity with
respect to $t$ is expressed in terms of a modulus of continuity.
Our goal is to find a sharp condition on the modulus of continuity
of the $\rho_{\alpha\beta}$'s, ensuring that $P$ enjoys the
backward uniqueness property in ${\mathcal H}_1^m$.

\par
Let $I\subset\R$ be a closed bounded interval, let ${\mathcal B}$ be a Banach space and let $f\colon I\to {\mathcal
B}$ be a continuous function. The {\it modulus of continuity of $f$} is the function $\mu(f,\cdot)\colon [0,1]\to
\R$ defined by
\begin{equation}\label{modcont1}
\mu(f,\tau):=\sup_{t,s\in I\atop 0\leq|t-s|\leq\tau}\|f(t)-f(s)\|_{\mathcal B}.
\end{equation}
Notice that $\mu(f,\cdot)$ is nondecreasing and $\mu(f,0)=0$. Since $f$
is uniformly continuous on $I$, it follows that
$\mu(f,\tau)\to0$ as $\tau\to 0$. If $f$ is nonconstant, then $\mu(f,\tau)>0$ for $\tau>0$ and $\mu(\cdot)=\mu(f,\cdot)$
satisfies
\begin{equation}\label{modcont2}
\sup_{t,s\in I\atop 0<|t-s|\leq1}{\frac{\|f(t)-f(s)\|_{\mathcal B}}{\mu(|t-s|)}}<+\infty.
\end{equation}
Moreover, $\mu(f,\cdot)$  is {\it minimal} with respect to this latter property, in the sense that, whenever a
nondecreasing function
$\mu\colon [0,1]\to\R$ satisfies (\ref{modcont2}), then $C\mu(\tau)\geq \mu(f,\tau)$ for some positive consant $C$.  It is
easy to check that
$\mu(f,\cdot)$ is sub-additive, that is
$\mu(f,\tau_1+\tau_2)\leq\mu(f,\tau_1)+\mu(f,\tau_2)$ whenever $0\leq\tau_1\leq\tau_2\leq\tau_1+\tau_2\leq1$. Then, by a
result of Efimov (\cite{Efi}, Lemma 4), there exists a concave (hence continuous) nondecreasing function
$\mu\colon[0,1]\to\R$ such that
$\mu(\tau)\leq\mu(f,\tau)\leq2\mu(\tau)$, $\tau\in[0,1]$. If $f$ is nonconstant, $\mu$
can be chosen to be strictly increasing. Therefore it is
natural to make the following

\begin{Def} Let $I\subset\R$ be a closed bounded interval, let ${\mathcal B}$ be a Banach space and
let $\mu\colon[0,1]\to[0,1]$ be a concave strictly increasing function, with $\mu(0)=0$. We say that a function $f\colon
I\to{\mathcal B}$ is
$\mu$-continuos (and we write $f\in {\mathcal C}^{\mu}(I,{\mathcal B})$) iff (\ref{modcont2}) is satisfied.
\end{Def}

Whenever $f\colon I\to {\mathcal B}$ is a continuous function, then certainly $f\in{\mathcal C}^\mu(I,{\mathcal B})$ for
some  concave strictly increasing function $\mu$, with $\mu(0)=0$ .
Moreover, if $f$ is nonconstant,
$\mu$ can be chosen in such a way that $\mu(\tau)\leq\mu(f,\tau)\leq2\mu(\tau)$, $\tau\in[0,1]$. If
$\mu(\tau)=\tau$, then ${\mathcal C}^\mu(I,{\mathcal B})={\rm Lip}(I,{\mathcal B})$; if $\mu(\tau)=\tau^\alpha$,
$0<\alpha<1$, then ${\mathcal C}^\mu(I,{\mathcal B})=C^\alpha(I,{\mathcal B})$.

\begin{Def}Let $\omega\colon[0,1]\to\R$ be an nondecreasing function, with $\omega(0)=0$. We say that $\omega$
satisfies the {\it Osgood condition} iff
\begin{equation}\label{Osgood}
\int_0^1{\frac1{\omega(s)}}\,ds=+\infty.
\end{equation}
\end{Def}
Condition (\ref{Osgood}) was introduced by Osgood in \cite{Osg}, while proving uniqueness for
ordinary differential equations with non-Lipschitz continuous nonlinearities. If $\omega(\tau)=\tau$, then $\omega$
satisfies the Osgood condition. If $\omega(\tau)=\tau^\alpha$, $0<\alpha<1$, then $\omega$ does not satisfy the Osgood
condition. If
$\omega(\tau)=\tau|\log\tau|$, then $\omega$ satisfies the Osgood condition.

The result of Mandache states that if a modulus of continuity $\omega$ does not satisfy the Osgood condition then there exists a parabolic operator of type (\ref{operator}) with $m=1$ such that the coefficients are
$C^\infty$ with respect to $x$ and $C^\omega $ with respect to $t$ and the backward uniqueness property does not hold.
\par Let us come
to our result. We shall consider here only
operators whose coefficients $\rho_{\alpha\beta}$'s are independent of the spatial variable $x$. If $m=1$, the general case
can be recovered by a microlocal approximation procedure similar to the one exploited in
\cite{CL}, provided the $\rho_{\alpha\beta}$'s are sufficiently smooth in $x$ (see
\cite{DSP} for details). We do not know whether it is possible to extend the result of \cite{DSP} to the case $m>1$.

\par If the coefficients are independent of $x$, the operator $P$ takes the simpler form
\begin{equation}\label{operator2}
P=\partial_t+{\hskip-4pt}\sum_{0\leq|\alpha|\leq 2m}{\hskip-4pt}
{\rm i}\,^{|\alpha|}\rho_{\alpha}(t)\partial_x^{\alpha},
\end{equation}
where $\rho_\alpha\in\R$ for all $\alpha$. Let $A$ be the set of all $n$-multiindices whose weight is smaller or equal than
$2m$, let $\ell$ be the cardinality of $A$ and let $R\colon[0,T]\to\R^\ell$, $R(t):=(\rho_\alpha(t))_{\alpha\in A}$, be a
continuous mapping.
Setting
\begin{equation}\label{H0}
\rho_k(t,\xi):=(-1)^k{{\sum_{|\alpha|=k}\rho_\alpha(t)\xi^\alpha}\over{|\xi|^k}},\quad
(t,\xi)\in[0,T]\times(\R^n\setminus\{0\}),\quad k=0,\dots,2m,
\end{equation}
we assume that the exists $\Lambda>0$ such that, for
$(t,\xi)\in[0,T]\times(\R^n\setminus\{0\})$,
\begin{equation}\label{H1}
|\rho_{k}(t,\xi)|\leq \Lambda,\quad k=0, \dots, 2m-1, \quad\text{and}\quad
1/\Lambda\leq
\rho_{2m}(t,\xi)\leq \Lambda.
\end{equation}

Consider the following {\it backward-parabolic} inequality:
\begin{equation}\label{bpi}
\|\partial_t u-{\hskip-4pt}\sum_{0\leq|\alpha|\leq 2m}{\hskip-4pt}
{\rm i}\,^{|\alpha|}\rho_{\alpha}(t)\partial_x^{\alpha}u\|_{L^2}\leq \tilde C\|u\|_{H^m}.
\end{equation}
The main result of the paper is the following

\begin{Theorem}\label{backuni}
Let the modulus of continuity $\mu(R,\cdot)$ of $R(\cdot):=(\rho_\alpha(\cdot))_{\alpha\in A}$ satisfy the Osgood
condition. If
$u\in{\mathcal H}_1^m$ satisfies (\ref{bpi}) and $u(0,\cdot)\equiv 0$ in $\R^n$, then $u\equiv 0$ in $[0,T]\times\R^n$.
\end{Theorem}

Notice that, if $R(\cdot)$ is constant, then $R\in{\mathcal C}^\mu([0,T], \R^\ell)$ with
$\mu(\tau)=\tau$. If
$R(\cdot)$ is nonconstant, then  we can find a concave strictly increasing function
$\mu\colon[0,1]\to\R$ such that
$\mu(\tau)\leq\mu(R,\tau)\leq2\mu(\tau)$, $\tau\in[0,1]$. It follows that in both cases $R\in{\mathcal
C}^\mu([0,T],\R^\ell)$ for some
$\mu$ which satisfies the Osgood condition. This observation is crucial for the proof of Theorem \ref{backuni}.

{\bf Remark 1.} Theorem \ref{backuni} allows to treat also operators with $x$-dependent coefficients up to the order $m$.
Indeed, all terms up to the order $m$ are absorbed by the right-hand side of the inequality (\ref{bpi}).

\par It is very likely that Theorem \ref{backuni} be sharp. Indeed not only the example of Mandache confirms it
in the case of $m=1$, but, by modifying a well known elliptic counterexample of
Pli\'s
\cite{Pl}, we can prove the following
\begin{Theorem}\label{counter} Let $\mu\colon[0,1]\to[0,1]$ be a concave strictly increasing function with $\mu(0)=0$. If
$\mu$ does not satisfy the Osgood  condition,
then for all $m\in\N$ there exist $l\in {\mathcal C}^\mu([0,1],\R)$, with $1/2\leq l(t)\leq 3/2$
for all
$t\in [0,1]$, $b_1$, $b_2$, $c\in C^\infty_b([0,1]
\times{\mathbb R}^2,\R)$, and  $u\in C^\infty_b([0,1]
\times{\mathbb R}^2,\R)$, with $u(1,\cdot)\equiv 0$ in
$\R^2$ but
$u\not\equiv 0$ in
$[0,1]\times\R^2$,  such that
\begin{multline}
\partial_t u+(-1)^m( \partial^{2m}_{x_1}u+ l(t)\partial^{2m}_{x_2}u)\\
+b_1(t,x)\partial_{x_1}u+b_2(t,x)\partial_{x_2}u+c(t,x)u=
0\quad \text{in $[0,1]\times{\mathbb R}^2$.}
\end{multline}
\end{Theorem}
{\bf Remark 2.} If $m=1$, we can take any function $\psi\in C^\infty(\R^2)$ such that $\psi(x)=e^{-|x|}$ for $|x|\geq1$,
and defining
$v(t,x):=\psi(x)u(t,x)$ we obtain a counterexample to the backward uniqueness property in ${\mathcal H}^1_1$. However if
$m>1$, by the same procedure we get only a non self-adjoint counterexample, with $x$-dependent coefficients up to the
order $2m-1$.

\par In the
next sections we  give sketches of the proofs of Theorems \ref{backuni} and  \ref{counter}.

\section {Proof of Theorem \ref{backuni}}
Theorem \ref{backuni} is a consequence of the following
\begin{Prop}\label{Carleman} Let $\mu\colon[0,1]\to[0,1]$ be a concave strictly increasing function with $\mu(0)=0$. Let
$T>0$ and let $R(\cdot)\in{\mathcal C}^\mu([0,T],\R^\ell)$, $R(\cdot):=(\rho_\alpha(\cdot))_{\alpha\in A}$, be a function
satisfying (\ref{H0})--(\ref{H1}). There exist
$C>0$,
$\gamma_0>0$ and a strictly increasing $C^2$-function $\Phi\colon[0,+\infty[\to[0,+\infty[$ such that, for all
$\gamma\geq\gamma_0$ and for all $u\in C^\infty_0(\R\times\R^n)$ with ${\rm supp}\,u\subset [0,T/2]\times\R^n$, the
following Carleman estimate holds:
\begin{multline}\label{carlest}
\int_0^{T/2}e^{{2\over\gamma}\Phi(\gamma(T-t))}\|\partial_t u-{\hskip-4pt}\sum_{0\leq|\alpha|\leq 2m}{\hskip-4pt}
{\rm i}\,^{|\alpha|}\rho_{\alpha}(t)\partial_x^{\alpha}u\|_{L^2}^2\,dt\\
\geq C \gamma^{1/2}\int_0^{T/2}e^{{2\over\gamma}\Phi(\gamma(T-t))}\|u\|^2_{H^m}\,dt.
\end{multline}
\end{Prop}

Let us briefly sketch how to prove Theorem \ref{backuni} from the Carleman estimate (\ref{carlest}).
First, we notice that if $\mu(R,\cdot)$ satisfies the Osgood condition then $R\in{\mathcal
C}^\mu([0,T],\R^\ell)$ for some concave strictly increasing function
$\mu$ which satisfies the Osgood condition. Second, by a density argument we have that (\ref{carlest}) holds
for any $u\in{\mathcal H}^m_1$ such that $u(0,\cdot)\equiv 0$ and $u(t,\cdot)\equiv0$ for $t\in[T/2,T]$. Now if $u\in
{\mathcal H}^m_1$ satisfies (\ref{bpi}) and $u(0,\cdot)\equiv 0$, we take $\vartheta\in C^\infty(\R)$, $\vartheta\equiv 0$
on $[T/2,+\infty]$, $\vartheta\equiv 1$ on $[0,T/3]$ and we apply (\ref{carlest}) to the function $\vartheta u$. We obtain
\begin{multline}\label{a1}
\int_0^{T/2}e^{{2\over\gamma}\Phi(\gamma(T-t))}\|\partial_t(\vartheta
u)-{\hskip-4pt}\sum_{0\leq|\alpha|\leq 2m}{\hskip-4pt}
{\rm i}\,^{|\alpha|}\rho_{\alpha}(t)\partial_x^{\alpha}
(\vartheta u)\|_{L^2}^2\,dt\\
\geq C \gamma^{1/2}\int_0^{T/2}e^{{2\over\gamma}\Phi(\gamma(T-t))}\|\vartheta u\|_{H^m}\,dt.
\end{multline}
Since $\vartheta\equiv 1$ for $t\in[0,T/3]$, (\ref{bpi}) and (\ref{a1}) imply
\begin{multline}\label{a2}
\int_{T/3}^{T/2}e^{{2\over\gamma}\Phi(\gamma(T-t))}\|\partial_t(\vartheta
u)-{\hskip-4pt}\sum_{0\leq|\alpha|\leq 2m}{\hskip-4pt}
{\rm i}\,^{|\alpha|}\rho_{\alpha}(t)\partial_x^{\alpha}
(\vartheta u)
\|_{L^2}^2\,dt\\
\geq(C\gamma^{1/2}-\tilde C)\int_0^{T/3}e^{{2\over\gamma}\Phi(\gamma(T-t))}\|u\|_{H^m}\,dt.
\end{multline}
Since $\Phi$ is increasing, for all sufficiently large $\gamma$ we have
\begin{equation}\label{a3}
\int_{T/3}^{T/2}\|\partial_t(\vartheta
u)-{\hskip-4pt}\sum_{0\leq|\alpha|\leq 2m}{\hskip-4pt}
{\rm i}\,^{|\alpha|}\rho_{\alpha}(t)\partial_x^{\alpha}
(\vartheta u)
\|_{L^2}^2\,dt\geq
{C\over 2}\gamma^{1/2}\int_0^{T/3}\| u\|^2_{H^m}\,dt.
\end{equation}
Letting $\gamma\to\infty$, we get $u\equiv 0$ in $[0,T/3]\times\R^n$. Finally, a standard connection argument implies
that $u\equiv 0$ in $[0,T]\times\R^n$.

\par Let us come to the proof of Lemma \ref{Carleman}. Let $\Phi\colon[0,+\infty[\to[0,+\infty[$ be of class $C^2$ and
increasing. Setting $v(t,x):=e^{{1\over\gamma}\Phi(\gamma(T-t))}u(t,x)$ and denoting by $\hat v(t,\xi)$ the Fourier
transform of $v(t,x)$ with respect to $x$, (\ref{carlest}) becomes
\begin{multline}\label{carlest2}
\int_0^{T/2}\int_{\R^n}|\partial_t\hat v(t,\xi)
-(\sum_{k=0}^{2m}\rho_k(t,\xi)|\xi|^k-\Phi'(\gamma(T-t)))\hat v(t,\xi)|^2\,\,d\xi \,dt\\
\geq C\gamma^{1/2}\int_0^{T/2}\int_{\R^n}(|\xi|^{2m}|\hat v(t,\xi)|^2+|\hat v(t,\xi)|^2)\,\,d\xi \,dt.
\end{multline}
Denoting by $\Xi$ the left member of (\ref{carlest2}), direct computation and integration by parts give
\begin{align*}
\Xi&=\int_0^{T/2}\int_{\R^n}|\partial_t\hat v(t,\xi)|^2\,\,d\xi \,dt\\
&+\int_0^{T/2}\int_{\R^n}| \sum_{k=0}^{2m}\rho_k(t,\xi)|\xi|^k    -\Phi'(\gamma(T-t))|^2|\hat v(t,\xi)|^2\,\,d\xi \,dt\\
&+\int_0^{T/2}\int_{\R^n}\gamma\Phi''(\gamma(T-t))|\hat v(t,\xi)|^2\,\,d\xi \,dt\\
&-2{\rm Re}\int_0^{T/2}\int_{\R^n}\partial_t\hat v(t,\xi) (\sum_{k=0}^{2m}\rho_k(t,\xi)|\xi|^k) \overline{\hat
v(t,\xi)}\,\,d\xi \,dt.
\end{align*}
If the $\rho_k(\cdot,\xi)$'s are Lipschitz continuos (that is: if $R(\cdot)$ is Lipschitz continuos),
one could just take $\Phi(\tau):=\tau^2$,
integrate by parts the double product and get the desired estimate (see \cite{LM} for details). If the
$\rho_k(\cdot,\xi)$'s are not Lipschitz continuos, we exploit a standard approximation procedure. We extend
$\rho_k(\cdot,\xi)$ on the whole
$\R$, we take
$\phi\in C^\infty_0(\R)$ such that
$\int_\R\phi(s)\,ds=1$,
$\phi\geq0$ and
${\rm supp}\,\phi\subset[-1/2,1/2]$, and then we define
\begin{equation}
\rho_{k,\epsilon}(t,\xi):=\int_\R \rho_k(s,\xi){1\over\epsilon}\phi\left({{t-s}\over{\epsilon}}\right)\,ds,\quad
(t,\xi)\in\R\times(\R^n\setminus\{0\}).
\end{equation}
It follows that $\rho_{k,\epsilon}(\cdot,\xi)\in C^\infty$ for every $\xi\in\R^n\setminus\{0\}$. Moreover
\begin{equation}
|\rho_{k,\epsilon}(t,\xi)-\rho_k(t,\xi)|\leq K\mu(\epsilon),\quad\text{$(t,\xi)\in\R\times(\R^n\setminus\{0\})$}
\end{equation}
and
\begin{equation}
|\rho'_{k,\epsilon}(t,\xi)|\leq K{{\mu(\epsilon)}\over\epsilon},\quad\text{$(t,\xi)\in\R\times(\R^n\setminus\{0\})$}
\end{equation}
(here ``$\,'\,$" indicates derivation with respect to $t$). Now let $\epsilon_1$, \dots, $\epsilon_{2m}$ be
approximation parameters to
be chosen later. Then, adding and subtracting
$\rho_{k,\epsilon_k}$ and integrating by parts with respect to $t$, we get:
\begin{multline}
-2{\rm Re}\int_0^{T/2}\int_{\R^n}\partial_t\hat v(t,\xi) (\sum_{k=0}^{2m}\rho_k(t,\xi)|\xi|^k)|\xi|^2\overline{\hat
v(t,\xi)}\,\,d\xi \,dt\\
\geq -\int_0^{T/2}\int_{\R^n}|\partial_t\hat v(t,\xi)|^2\,\,d\xi \,dt -K\int_0^{T/2}\int_{\R^n}
(\sum_{k=0}^{2m}{{\mu(\epsilon_k)}\over{\epsilon_k}}|\xi|^k)
|\hat v(t,\xi)|^2\,\,d\xi \,dt\\ -K^2\int_0^{T/2}\int_{\R^n} (\sum_{k=0}^{2m}\mu(\epsilon_k)^2|\xi|^{2k})
|\hat v(t,\xi)|^2\,\,d\xi \,dt\\
\end{multline}
Now the first key idea is to let the approximation parameters $\epsilon_k$ depend on $\xi$ (cf \cite{CDGS}).
First, we observe that, by (\ref{H1}), there exist $N_0\geq 1$ and $\Lambda_0>0$ such that, for all $|\xi|\geq N_0$,
\begin{equation}\label{H2}
{{1}\over{\Lambda_0}}|\xi|^{2m}\leq \sum_{k=0}^{2m}\rho_k(t,\xi)|\xi|^k \leq\Lambda_0|\xi|^{2m}
\end{equation}
Then we take
\begin{equation}
\epsilon_k:=\begin{cases}|\xi|^{-k}&\text{if $|\xi|\geq N_0$}\\
N_0^{-k}&\text{if $|\xi|\leq N_0$}
\end{cases}
\end{equation}
With this choice, noticing also that $s^2\mu(1/s)$ is increasing on $[1,+\infty]$, we obtain that there exists a positive
constant $\tilde K$ such that:
\begin{align*}
&\int_0^{T/2}\int_{\R^n}|\partial_t\hat v(t,\xi)
-(\sum_{k=0}^{2m}\rho_k(t,\xi)|\xi|^k-\Phi'(\gamma(T-t)))\hat v(t,\xi)|^2\,\,d\xi \,dt&&\\
&\geq \gamma\int_0^{T/2}\int_{\R^n}\Phi''(\gamma(T-t))|\hat v(t,\xi)|^2\,\,d\xi \,dt&&(i)\\
&+\int_0^{T/2}\int_{\R^n}| \sum_{k=0}^{2m}\rho_k(t,\xi)|\xi|^k-\Phi'(\gamma(T-t))|^2|\hat v(t,\xi)|^2\,\,d\xi \,dt&&(ii)\\
&-\tilde K \int_0^{T/2}\int_{\R^n}|\hat v(t,\xi)|^2\,\,d\xi\,dt&&(iii)\\
&-\tilde K\int_0^{T/2}\int_{\{|\xi|\geq N_0\}}\mu(1/|\xi|^{2m})|\xi|^{4m}|\hat v(t,\xi)|^2\,\,d\xi \,dt&&(iv)\\
\end{align*}
Now we observe that:
\begin{itemize}
\item the summand $(i)$ behaves well, provided $\Phi''(\tau)\geq1$ for large $\tau$;
\item the same condition let the summand $(iii)$ be absorbed by $(i)$;
\item when $|\xi|^{2m}\geq 2\Lambda_0\Phi'(\gamma(T-t)))$, then the integrand in $(ii)$ behaves like
$|\xi|^{4m}$, which is enough to compensate the integrand in $(iv)$;
\item when  $N_0^{2m}\leq |\xi|^{2m}\leq(1/2\Lambda_0) \Phi'(\gamma(T-t)))$, then again the integrand in $(ii)$ behaves
like
$|\xi|^{4m}$, which is enough to compensate the integrand in $(iv)$;
\item the difficult case is when $|\xi|^{2m} \sim\Phi'(\gamma(T-t))$.
\end{itemize}

At this point the second key idea is to modulate the weight $\Phi$ on the function $\mu$ (cf Tarama \cite{tara}). Roughly
speaking, we ask that, when
$|\xi|^{2m} \sim\Phi'(\gamma(T-t))$, then the integrand in $(iv)$ must be compensated by the integrand in $(i)$. More
precisely, we ask that $\Phi''(\gamma(T-t))\sim\mu(1/|\xi|^{2m})|\xi|^{4m} $. In other words, the Carleman estimates
(\ref{carlest}) will follow, provided $\Phi$ satisfies the ordinary differential equation
\begin{equation}\label{ode}
\Phi''=\mu(1/\Phi')(\Phi')^2.
\end{equation}
All we have to do then is to find a solution of (\ref{ode}) and to check that:
\begin{itemize}
\item $\Phi$ is defined on $[0+\infty[$, i.e. it does not blow up in finite time;
\item $\Phi$ is positive and increasing;
\item $\Phi''(\tau)\geq 1$ for all sufficiently large $\tau$.
\end{itemize}
Equation (\ref{ode}) can be easily solved by separation of variables. The explicit solution
of the Cauchy problem with initial values $\Phi(0)=0$ and $\Phi'(0)=1$ is given by:
\begin{align*}
\eta(t)&:=\int^1_{1/t}{{1}\over{\mu(s)}}\,ds,&&t\geq 1\\
\Phi(\tau)&:=\int_0^\tau\eta^{-1}(r)dr,&&\tau\geq 0
\end{align*}
The {\it Osgood condition} precisely guarantees that $\Phi$ is defined on $[0,+\infty[$. The other properties that
we require for $\Phi$ follow by easy computation. With this choice of $\Phi$ we finally get the desired Carleman estimate
(\ref{carlest}). The details are left to the reader.

\section {Proof of Theorem \ref{counter}}

The proof of Theorem 2 is very similar to that one of Theorem 3 in \cite{DSP}. Also in this case we will follow closely the construction of the example in \cite {Pl}. Let $A$, $B$, $C$, $J$
be four $C^\infty$ functions defined in $\mathbb R$ with $0\leq
A(s),\ B(s),\ C(s)\leq 1$, $-2\leq J(s)\leq 2$ for all $s\in
{\mathbb R}$ and $$
\begin{array}{ll}
\displaystyle{A(s)=1\quad{\rm for}\ s\leq {1\over 5}, }&
\quad\displaystyle{A(s)=0\quad{\rm for}\
s\geq {1\over 4},}\\[0.3 cm]
\displaystyle{B(s)=0\quad{\rm for}\ s\leq 0\ {\rm or}\ s\geq 1, }&\quad
\displaystyle{B(s)=1\quad{\rm for}\ {1\over 6}\leq s\leq {1\over 2},}\\[0.3 cm]
\displaystyle{C(s)=0\quad{\rm for}\ s\leq {1\over 4}, }&
\quad\displaystyle{C(s)=1\quad{\rm for}\
s\geq {1\over 3},}\\[0.3 cm]
\displaystyle{J(s)=-2\quad{\rm for}\ s\leq {1\over 6}\ {\rm or}\ s\geq
{1\over 2}, }&\quad
\displaystyle{J(s)=2\quad{\rm for}\ {1\over 5}\leq s\leq {1\over 3}.}\\
\end{array}
$$ Let $(a_n)_n$, $(z_n)_n$ be two real sequences such that
\begin{eqnarray}
\displaystyle{-1<a_n<a_{n+1}\quad {\rm for\ all}\ n\geq
1},&&\displaystyle{\lim_n a_n=0,}\label
{4.1}\\[0.3 cm]
\displaystyle{1<z_n<z_{n+1}\quad {\rm for\ all}\ n\geq
1},&&\displaystyle{\lim_n
z_n=+\infty;}\label {4.2}
\end{eqnarray} and let us define $r_n=a_{n+1}-a_n$, $q_1=0$,
$q_n=\sum_{k=2}^nz_k r_{k-1}$ for all
$n\geq 2$, and $p_n=(z_{n+1}-z_n)r_n$. We suppose moreover that
\begin{equation} p_n>1\quad {\rm for \ Êall}\ n\geq 1.
\label{4.3}
\end{equation} We set $A_n(t)=A({t-1-a_n\over r_n})$, $B_n(t)=B({t-1-a_n\over
r_n})$,
$C_n(t)=C({t-1-a_n\over r_n})$ and $J_n(t)=J({t-1-a_n\over r_n})$. We define
$$
\begin{array}{l}
\displaystyle{v_n(t,x_1)=\exp (-q_n-z_n(t-1-a_n))\cos \root {2m}\of {z_n} x_1,}\\[0.3 cm]
\displaystyle{w_n(t,x_2)=\exp (-q_n-z_n(t-1-a_n)+J_n(t)p_n)\cos\root {2m}\of{z_n}
x_2,}\\
\end{array}
$$ and
$$
\begin{array}{l} u(t,x_1, x_2)\\[0.3 cm] =\left\{
\begin{array}{ll}
\displaystyle{v_1(t, x_1)} &\quad\displaystyle{{\rm for }\ 0\leq t\leq 1+a_1, }\\[0.3 cm]
\displaystyle{A_n(t)v_n(t,x_1)+B_n(t)w_n(t,x_2)}&\\[0.3 cm]
\qquad\qquad\qquad\qquad\displaystyle{+C_n(t)v_{n+1}(t,x_1)}&\quad\displaystyle{{\rm
for }\
1+a_n\leq t\leq 1+a_{n+1} ,}\\[0.3 cm]
\displaystyle{0} &\quad\displaystyle{{\rm for }\ t=1. }\\
\end{array}
\right.
\end{array}
$$ If for all $\alpha$, $\beta$ $\gamma>0$
\begin{equation}
\lim_n\exp(-q_n+2p_n)z_{n+1}^\alpha p_n^\beta r_n^{-\gamma}=0
\label{4.4}
\end{equation} then $u$ is a $C^\infty_b([0,1]\times{\mathbb R}^2, {\mathbb R})$ function. We define
$$ l(t)=\left\{
\begin{array}{ll}
\displaystyle{1} &\displaystyle{{\rm for }\ t\leq 1+a_1\ {\rm or }\ t=1,
}\\[0.3 cm]
\displaystyle{1+J'_n(t)p_nz^{-1}_n}&\displaystyle{{\rm for }\ 1+a_n\leq t\leq
1+a_{n+1} .}\\
\end{array}
\right.
$$ The condition
\begin{equation}
\sup_n\; \{p_nr_n^{-1}z_n^{-1}\}\leq {1\over 2\|J'\|_{L^\infty}}
\label{4.5}
\end{equation} guarantees that the operator $\text{\it
\L}=\partial_t+(-1)^{m}(\partial^{2m}_{x_1}-l(t)\partial^{2m}_{x_2})$ is parabolic.
Moreover $l$ is a $C^\mu$
function under the condition
\begin{equation}
\sup_n\; \{{p_nr_n^{-1}z_n^{-1}\over \mu(r_n)}\}<+\infty.
\label{4.6}
\end{equation} Finally we define
$$
\begin{array}{l}
\displaystyle{b_1=- {\text{\it \L} u\over
u^2+(\partial_{x_1}u)^2+(\partial_{x_2}u)^2}\partial_{x_1}u,}\\[0.3cm]
\displaystyle{b_2=- {\text{\it \L} u\over
u^2+(\partial_{x_1}u)^2+(\partial_{x_2}u)^2}\partial_{x_2}u,}\\[0.3cm]
\displaystyle{c=- {\text{\it \L} u\over
u^2+(\partial_{x_1}u)^2+(\partial_{x_2}u)^2}u}.\\
\end{array}
$$
As in \cite {Pl},  or similarly in \cite{DSP}, the coefficients $b_1$, Ê$b_2$, Ê$c$ will be in
$C^\infty_b$ if for all
$\alpha$, $\beta$, $\gamma>0$
\begin{equation}
\lim_n\exp(-p_n)z_{n+1}^\alpha p_n^\beta r_n^{-\gamma}=0.
\label{4.7}
\end{equation} We choose
\begin{equation} a_n= -\sum _{j=n}^{+\infty}{1\over (j+k_0)^2\mu({1\over
j+k_0})}, \qquad
z_n=(n+k_0)^3
\label{4.8}
\end{equation} with $k_0$ sufficiently large.

To conclude the proof it will be sufficient to verify in the same way than in   \cite{DSP} that, with the choice (\ref{4.8}), the conditions
(\ref{4.1}),...,
(\ref{4.7}) hold. We let it to the reader.

\end{document}